\documentclass[12pt]{article}

\usepackage{amsmath,amssymb,amsthm,amscd}

\begin{document}

\newcommand{\End}{{\rm{End}\ts}}
\newcommand{\Hom}{{\rm{Hom}}}
\newcommand{\ch}{{\rm{ch}\ts}}
\newcommand{\non}{\nonumber}
\newcommand{\wt}{\widetilde}
\newcommand{\wh}{\widehat}
\newcommand{\ot}{\otimes}
\newcommand{\la}{\lambda}
\newcommand{\La}{\Lambda}
\newcommand{\al}{\alpha}
\newcommand{\be}{\beta}
\newcommand{\ga}{\gamma}
\newcommand{\si}{\sigma}
\newcommand{\vp}{\varphi}
\newcommand{\de}{\delta^{}}
\newcommand{\om}{\omega^{}}
\newcommand{\hra}{\hookrightarrow}
\newcommand{\ve}{\varepsilon}
\newcommand{\ts}{\,}
\newcommand{\qin}{q^{-1}}
\newcommand{\tss}{\hspace{1pt}}
\newcommand{\U}{ {\rm U}}
\newcommand{\Y}{ {\rm Y}}
\newcommand{\CC}{\mathbb{C}\tss}
\newcommand{\ZZ}{\mathbb{Z}\tss}
\newcommand{\Z}{\mathbb{Z}}
\newcommand{\A}{\mathcal{A}}
\newcommand{\Pc}{\mathcal{P}}
\newcommand{\Qc}{\mathcal{Q}}
\newcommand{\Tc}{\mathcal{T}}
\newcommand{\Bc}{\mathcal{B}}
\newcommand{\Ec}{\mathcal{E}}
\newcommand{\Hc}{\mathcal{H}}
\newcommand{\Ar}{{\rm A}}
\newcommand{\Ir}{{\rm I}}
\newcommand{\Zr}{{\rm Z}}
\newcommand{\gl}{\mathfrak{gl}}
\newcommand{\Pf}{{\rm Pf}}
\newcommand{\oa}{\mathfrak{o}}
\newcommand{\spa}{\mathfrak{sp}}
\newcommand{\g}{\mathfrak{g}}
\newcommand{\ka}{\mathfrak{k}}
\newcommand{\p}{\mathfrak{p}}
\newcommand{\sll}{\mathfrak{sl}}
\newcommand{\agot}{\mathfrak{a}}
\newcommand{\qdet}{ {\rm qdet}\ts}
\newcommand{\sdet}{ {\rm sdet}\ts}
\newcommand{\sgn}{ {\rm sgn}\ts}
\newcommand{\Sym}{\mathfrak S}
\newcommand{\fand}{\quad\text{and}\quad}
\newcommand{\Fand}{\qquad\text{and}\qquad}

\renewcommand{\theequation}{\arabic{section}.\arabic{equation}}

\newtheorem{thm}{Theorem}[section]
\newtheorem{lem}[thm]{Lemma}
\newtheorem{prop}[thm]{Proposition}
\newtheorem{cor}[thm]{Corollary}
\newtheorem{conj}[thm]{Conjecture}

\theoremstyle{definition}
\newtheorem{defin}[thm]{Definition}
\newtheorem{example}[thm]{Example}

\theoremstyle{remark}
\newtheorem{remark}[thm]{Remark}

\newcommand{\bth}{\begin{thm}}
\renewcommand{\eth}{\end{thm}}
\newcommand{\bpr}{\begin{prop}}
\newcommand{\epr}{\end{prop}}
\newcommand{\ble}{\begin{lem}}
\newcommand{\ele}{\end{lem}}
\newcommand{\bco}{\begin{cor}}
\newcommand{\eco}{\end{cor}}
\newcommand{\bde}{\begin{defin}}
\newcommand{\ede}{\end{defin}}
\newcommand{\bex}{\begin{example}}
\newcommand{\eex}{\end{example}}
\newcommand{\bre}{\begin{remark}}
\newcommand{\ere}{\end{remark}}
\newcommand{\bcj}{\begin{conj}}
\newcommand{\ecj}{\end{conj}}

\newcommand{\bal}{\begin{aligned}}
\newcommand{\eal}{\end{aligned}}
\newcommand{\beq}{\begin{equation}}
\newcommand{\eeq}{\end{equation}}
\newcommand{\ben}{\begin{equation*}}
\newcommand{\een}{\end{equation*}}

\newcommand{\bpf}{\begin{proof}}
\newcommand{\epf}{\end{proof}}

\def\beql#1{\begin{equation}\label{#1}}

\title{\Large\bf Symmetries and invariants of twisted quantum
algebras and associated Poisson algebras}

\author{{A. I. Molev\quad and\quad E. Ragoucy}}

\date{} 
\maketitle

\vspace{5 mm}

\begin{abstract}
We construct an action of the braid group $B_N$
on the twisted quantized enveloping algebra $\U'_q(\oa_N)$
where the elements of $B_N$ act as automorphisms.
In the classical limit
$q\to 1$ we recover the action of $B_N$ on the
polynomial functions on the space of upper triangular matrices with
ones on the diagonal.
The action preserves the Poisson bracket
on the space of polynomials which was introduced by Nelson and Regge
in their study of quantum gravity and re-discovered
in the mathematical literature.
Furthermore, we construct a Poisson bracket on the
space of polynomials associated with another
twisted quantized enveloping algebra $\U'_q(\spa_{2n})$.
We use the Casimir elements of both
twisted quantized enveloping algebras
to re-produce some well-known and construct some new
polynomial invariants
of the corresponding Poisson algebras.

\vspace{5 mm}

Preprint LAPTH-1174/07
\end{abstract}


\vspace{30 mm}

\noindent
School of Mathematics and Statistics\newline
University of Sydney,
NSW 2006, Australia\newline
alexm@maths.usyd.edu.au

\vspace{7 mm}

\noindent
LAPTH, Chemin de Bellevue, BP 110\newline
F-74941 Annecy-le-Vieux cedex, France\newline
ragoucy@lapp.in2p3.fr

\newpage

\section{Introduction}\label{sec:int}
\setcounter{equation}{0}

Deformations of
the commutation relations of the orthogonal Lie algebra $\oa_3$
were considered by many authors.
The earliest reference we are aware of is Santilli~\cite{s:ru}. Such deformed
relations can be written as
\beql{ortdef}
q\tss XY-YX=Z,\qquad
q\tss YZ-ZY=X,\qquad
q\tss ZX-XZ=Y.
\eeq
More precisely, regarding $q$ as a formal variable,
we consider the associative algebra $\U'_q(\oa_3)$ over the field
of rational functions $\CC(q)$ in $q$ with the generators $X,Y,Z$ and
defining relations \eqref{ortdef}.
From an alternative viewpoint, relations \eqref{ortdef} define
a family of algebras depending on the complex parameter $q$.
The same algebras were also defined by Odesskii~\cite{o:as},
Fairlie~\cite{f:qd} and Nelson, Regge and Zertuche~\cite{nrz:hg}.
Putting $q=1$ in \eqref{ortdef} we get
the defining relations of the universal
enveloping algebra $\U(\oa_3)$. The algebra $\U'_q(\oa_3)$ should be
distinguished from the quantized enveloping algebra
$\U_q(\oa_3)\cong\U_q(\sll_2)$. The latter is
a deformation of $\U(\oa_3)$ in the class of Hopf algebras;
see e.g. Chari and Pressley~\cite[Section~6]{cp:gq}.

Introducing the generators
\ben
x=(q-\qin)\tss{X},\qquad y=(q-\qin)\tss{Y},\qquad z=(q-\qin)\tss{Z},
\een
we can write the defining relations of $\U'_q(\oa_3)$ in the
equivalent form
\ben
\bal
q\tss xy-yx&=(q-\qin)\tss  z,\\
q\tss yz-zy&=(q-\qin)\tss  x,\\
q\tss zx-xz&=(q-\qin)\tss  y.
\eal
\een
Note that the element
$
x^2+q^{-2}\ts y^2+z^2-xyz
$
belongs to the center of $\U'_q(\oa_3)$.
This time, putting $q=1$ into the defining relations we get
the algebra of polynomials $\CC[x,y,z]$.
Moreover, this algebra can be equipped with a Poisson bracket
in a usual way
\ben
\{f,g\}=\frac{fg-gf}{1-q}\ts\Big|_{q=1}.
\een
Thus, $\CC[x,y,z]$ becomes a {\it Poisson algebra\/}
with the bracket given by
\beql{dubropb}
\{x,y\}=xy-2z,\qquad
\{y,z\}=yz-2x,\qquad
\{z,x\}=zx-2y.
\eeq
These formulas are contained in the paper
by Nelson, Regge and Zertuche~\cite{nrz:hg}.
In the classical limit $q\to 1$
the central element $x^2+q^{-2}\ts y^2+z^2-xyz$
becomes the {\it Markov polynomial}
$x^2+y^2+z^2-xyz$ which is
an invariant of the bracket.
The Poisson bracket \eqref{dubropb}
was re-discovered by Dubrovin~\cite{d:gt}, where
$x,y,z$ are interpreted as the entries of
$3\times 3$ upper triangular matrices with
ones on the diagonal
(the {\it Stokes matrices\/})
\ben
\begin{pmatrix} 1&x&y\\
                0&1&z\\
                0&0&1
\end{pmatrix}.
\een

For an arbitrary $N$ the {\it twisted quantized enveloping algebra\/}
$\U'_q(\oa_N)$ was introduced by Gavrilik and Klimyk~\cite{gk:qd}
which essentially coincides with the algebra of
Nelson and Regge~\cite{nr:qg}.
Both in the orthogonal and symplectic case the
twisted analogues of the quantized enveloping algebras were
introduced by Noumi~\cite{n:ms} using
an $R$-matrix approach. In the orthogonal case
this provides an alternative presentation of $\U'_q(\oa_N)$.
The finite-dimensional irreducible representations
of the algebra $\U'_q(\oa_N)$ were classified
by Iorgov and Klimyk~\cite{ik:ct}.

In the limit $q\to 1$
the twisted quantized enveloping algebra $\U'_q(\oa_N)$ gives rise to
a Poisson algebra of polynomial functions $\Pc_N$ on the space of Stokes
matrices. The corresponding Poisson bracket
was given in \cite{nr:qg}.
The same bracket was also found by Ugaglia~\cite{u:ps}, Boalch~\cite{b:sm}
and Bondal~\cite{b:sg, b:sgr}.  This Poisson structure was studied
by Ping Xu~\cite{x:ds}
in the context of Dirac submanifolds, while
Chekhov and Fock~\cite{cf:od} considered it in relation with
the Teichm\"{u}ller spaces.
A quantization of the Poisson algebra of Stokes matrices
leading to the algebra $\U'_q(\oa_N)$
was constructed by Ciccoli and Gavarini~\cite{cg:qd}
in the context of the general ``quantum quality principle";
see also Gavarini~\cite{g:pb}.
It was shown by
Odesskii and Rubtsov~\cite{or:pp} that the Poisson bracket
on the space of Stokes matrices
is essentially determined
by its Casimir elements.

Automorphisms of both the algebra $\U'_q(\oa_N)$
and the Poisson bracket on $\Pc_N$
were given in \cite{nr:gg, nr:ig}, although
the explicit group relations between them were only
discussed in the classical limit for $N=6$.
An action of the braid group $B_N$ on the Poisson algebra $\Pc_N$
was given by Dubrovin~\cite{d:gt} and Bondal~\cite{b:sg}.

In this paper we produce a ``quantized" action
of $B_N$ on the twisted quantized enveloping algebra $\U'_q(\oa_N)$,
where the elements of $B_N$ act
as automorphisms.
Since $\U'_q(\oa_N)$
is a subalgebra of the quantized enveloping algebra $\U_q(\gl_N)$,
one could expect that Lusztig's action of
$B_N$ on $\U_q(\gl_N)$ (see \cite{l:fd}) leaves the subalgebra $\U'_q(\oa_N)$
invariant. However, this turns out not to be true, and the
action of $B_N$ on $\U'_q(\oa_N)$ can rather be regarded
as a $q$-version of the natural action of the symmetric group $\Sym_N$
on the universal enveloping algebra $\U(\oa_N)$.

The relationship between $\U'_q(\oa_N)$
and the Poisson algebra $\Pc_N$
can also be exploited in a different way.
Some families of Casimir elements of $\U'_q(\oa_N)$ were produced
by Noumi, Umeda and Wakayama~\cite{nuw:dp},
Gavrilik and Iorgov~\cite{gi:ce} and
Molev, Ragoucy and Sorba~\cite{mrs:cs}. This gives
the respective families of Casimir elements of the
Poisson algebra.
We show that
the Casimir elements of \cite{mrs:cs} specialize precisely to
the coefficients of the characteristic polynomial
of Nelson and Regge~\cite{nr:ig}. This polynomial
was re-discovered by Bondal~\cite{b:sg} who also produced
an algebraically
independent set of generators of the
subalgebra of invariants of the Poisson algebra $\Pc_N$.
Furthermore, using \cite{gi:ce} and \cite{nuw:dp} we obtain new Pfaffian
type invariants and analogues of the Gelfand invariants.

In a similar manner, we use the twisted quantized enveloping algebra
$\U'_q(\spa_{2n})$ associated with the symplectic Lie algebra $\spa_{2n}$
to produce a symplectic version of the above results.
First, we construct a Poisson algebra associated with $\U'_q(\spa_{2n})$
by taking the limit $q\to 1$ and thus produce explicit
formulas for the Poisson bracket on
the corresponding space of matrices. Then using the Casimir
elements of $\U'_q(\spa_{2n})$ constructed in \cite{mrs:cs},
we produce a family of invariants of the Poisson algebra
analogous to \cite{b:sg} and \cite{nr:ig}.
We also show that some elements of the braid group $B_{2n}$
preserve the subalgebra $\U'_q(\spa_{2n})$ of $\U_q(\gl_{2n})$.
We conjecture that there exists an action
of the semi-direct product $B_n\ltimes\ZZ^n$
on $\U'_q(\spa_{2n})$ analogous to the $B_N$-action on $\U'_q(\oa_{N})$.
We show that the conjecture is true for $n=2$.

This work was inspired by
Alexei Bondal's talk at the Prague's conference ISQS 2006.
We would like to thank Alexei
for many stimulating discussions.
The financial support of the Australian
Research Council is acknowledged. The second author
is grateful to the University of Sydney for the warm hospitality
during his visit.

After we prepared the first version
of our paper we learned of a recent preprint by L.~Chekhov~\cite{ch:tt}
where he produces (without detailed proofs)
an action of the braid group $B_N$ on $\U'_q(\oa_{N})$ equivalent to ours.

\section{Braid group action}\label{sec:bga}
\setcounter{equation}{0}

We start with some definitions and recall some well-known results.
Let $q$ be a formal variable.
The {\it quantized enveloping algebra\/} $\U_q(\gl_N)$ is
an algebra over $\CC(q)$ generated
by elements $t_{ij}$ and $\bar t_{ij}$ with $1\leqslant i,j\leqslant N$
subject to the relations
\beql{defrel}
\bal
t_{ij}&=\bar t_{ji}=0, \qquad 1 \leqslant i<j\leqslant N,\\
t_{ii}\ts \bar t_{ii}&=\bar t_{ii}\ts t_{ii}=1,\qquad 1\leqslant i\leqslant N,\\
R\ts T_1T_2&=T_2T_1R,\qquad R\ts \overline T_1\overline T_2=
\overline T_2\overline T_1R,\qquad
R\ts \overline T_1T_2=T_2\overline T_1R.
\eal
\end{equation}
Here $T$ and $\overline T$ are the matrices
\beql{matrt}
T=\sum_{i,j}t_{ij}\ot E_{ij},\qquad \overline T=\sum_{i,j}
\overline t_{ij}\ot E_{ij},
\end{equation}
which are regarded as elements of the algebra $\U_q(\gl_N)\ot \End\CC^N$,
the $E_{ij}$ denote the standard matrix units and the indices run over
the set $\{1,\dots,N\}$.
Both sides of each of the $R$-matrix relations in \eqref{defrel}
are elements of $\U_q(\gl_N)\ot \End\CC^N\ot \End\CC^N$ and the subscripts
of $T$ and $\overline T$ indicate the copies of $\End\CC^N$,
e.g.,
\ben
T_1=\sum_{i,j}t_{ij}\ot E_{ij}\ot 1,\qquad
T_2=\sum_{i,j}t_{ij}\ot 1\ot E_{ij},
\een
while $R$ is
the $R$-matrix
\beql{rmatrixc}
R=q\ts\sum_i E_{ii}\ot E_{ii}+\sum_{i\ne j} E_{ii}\ot E_{jj}+
(q-\qin)\sum_{i<j}E_{ij}\ot E_{ji}.
\end{equation}

In terms of the
generators the defining relations between the $t_{ij}$
can be written as
\beql{defrelg}
q^{\delta_{ij}}\ts t_{ia}\ts t_{jb}-
q^{\delta_{ab}}\ts t_{jb}\ts t_{ia}
=(q-\qin)\ts (\de_{b<a} -\de_{i<j})
\ts t_{ja}\ts t_{ib},
\end{equation}
where $\de_{i<j}$ equals $1$ if $i<j$, and $0$ otherwise.
The relations between the $\bar t_{ij}$
are obtained by replacing $t_{ij}$ by $\bar t_{ij}$ everywhere in
\eqref{defrelg}, while the relations involving both
$t_{ij}$ and $\bar t_{ij}$ have the form
\beql{defrelg2}
q^{\delta_{ij}}\ts \bar t_{ia}\ts t_{jb}-
q^{\delta_{ab}}\ts t_{jb}\ts \bar t_{ia}
=(q-\qin)\ts (\de_{b<a}\ts t_{ja}\ts \bar t_{ib} -\de_{i<j}\ts
\ts \bar t_{ja}\ts t_{ib}).
\end{equation}

The {\it braid group\/} $B_N$ is generated
by elements $\be_1,\dots,\be_{N-1}$ subject to the defining
relations
\ben
\be_i\tss \be_{i+1}\tss \be_i=\be_{i+1}\tss \be_i\tss \be_{i+1},
\qquad i=1,\dots,N-2
\een
and
\ben
\be_i\tss \be_j=\be_j\tss \be_i,\qquad |i-j|>1.
\een
The group $B_N$ acts on the algebra $\U_q(\gl_N)$
by automorphisms; see Lusztig~\cite{l:fd}. Explicit formulas
for the images of the generators are found from \cite{l:fd}
by re-writing the action in terms of the presentation
\eqref{defrel}. For any $i=1,\dots,N-1$ we have
\ben
\be_i:
    t_{ii}\mapsto t_{i+1,i+1},\qquad
    t_{i+1,i+1}\mapsto t_{ii},\qquad
    t_{kk}\mapsto t_{kk} \qquad\text{if}\quad k\ne i,i+1,
\een
\ben
\bal
\be_i:{}&{}t_{i+1,i}\mapsto \qin\ts \bar t_{i,i+1}\ts t_{ii}^2\\
    &t_{ik}\mapsto q\ts t_{ik}\ts t_{i+1,i}\ts \bar t_{ii}
    -t_{i+1,k},
    \qquad
    &&t_{i+1,k}\mapsto \qin\ts t_{ik},
    \qquad&&\text{if}\quad k\leqslant i-1\\
    &t_{li}\mapsto \qin\ts \bar t_{i,i+1}\ts t_{li}\ts t_{ii}-t_{l,i+1},
    \qquad
    &&t_{l,i+1}\mapsto q\ts t_{li},
    \qquad&&\text{if}\quad l\geqslant i+2\\
    &t_{kl}\mapsto t_{kl}
    &&\text{in all remaining cases,}
\eal
\een
and
\ben
\bal
\be_i:{}&{}\bar t_{i,i+1}\mapsto q\ts \bar t_{ii}^{\ts\ts 2} \ts t_{i+1,i}\\
    &\bar t_{ki}\mapsto \qin\ts t_{ii}\ts \bar t_{i,i+1}\bar t_{ki}
    -\bar t_{k,i+1},
    \qquad
    &&\bar t_{k,i+1}\mapsto q\ts \bar t_{ki},
    \qquad&&\text{if}\quad k\leqslant i-1\\
    &\bar t_{il}\mapsto q\ts \bar t_{ii}\ts \bar t_{il}\ts t_{i+1,i}
    -\bar t_{i+1,l},
    \qquad
    &&\bar t_{i+1,l}\mapsto \qin\ts \bar t_{il},
    \qquad&&\text{if}\quad l\geqslant i+2\\
    &\bar t_{kl}\mapsto \bar t_{kl}
    &&\text{in all remaining cases.}
\eal
\een

Following Noumi~\cite{n:ms} we define the {\it twisted
quantized enveloping algebra\/} $\U'_q(\oa_N)$
as the subalgebra of $\U_q(\gl_N)$ generated by the matrix
elements $s_{ij}$ of the matrix $S=T\ts \overline T^{\ts t}$ so that
\ben
s_{ij}=\sum_{k=1}^N t_{ik}\ts\bar t_{jk}.
\een
Equivalently, $\U'_q(\oa_N)$ is generated by the
elements $s_{ij}$ subject only to the relations
\begin{align}\label{sijo}
s_{ij}&=0, \qquad 1 \leqslant i<j\leqslant N,\\
\label{sii1}
s_{ii}&=1,\qquad 1\leqslant i\leqslant N,\\
\label{rsrs}
R\ts S_1& R^{\ts t} S_2=S_2R^{\ts t} S_1R,
\end{align}
where $R^{\ts t}:=R^{\ts t_1}$ denotes the element obtained from $R$ by
the transposition in the first tensor factor:
\beql{rt}
R^{\ts t}=q\ts\sum_i E_{ii}\ot E_{ii}+\sum_{i\ne j} E_{ii}\ot E_{jj}+
(q-\qin)\sum_{i<j}E_{ji}\ot E_{ji}.
\end{equation}
In terms of the generators, the relations \eqref{rsrs} take the form
\beql{drabs}
\bal
q^{\delta_{jk}+\delta_{ik}}\ts s_{ij}\ts s_{kl}-
q^{\delta_{jl}+\delta_{il}}\ts s_{kl}\ts s_{ij}
{}&=(q-\qin)\ts q^{\delta_{ji}}\ts (\de_{l<j} -\de_{i<k})
\ts s_{kj}\ts s_{il}\\
{}&+(q-\qin)\ts \big(q^{\delta_{jl}}\ts \de_{l<i}\ts s_{ki}\ts s_{lj}
- q^{\delta_{ik}}\ts \de_{j<k}\ts s_{ik}\ts s_{jl}\big)\\
{}&+ (q-\qin)^2\ts  (\de_{l<j<i} -\de_{j<i<k})\ts s_{ki}\ts s_{jl},
\eal
\end{equation}
where $\de_{i<j}$ or $\de_{i<j<k}$ equals $1$ if the subscript inequality
is satisfied, and $0$ otherwise.
Equivalently, the set of relations can also be written as
\beql{drabssi}
\bal
s_{ij}\ts s_{kl}-s_{kl}\ts s_{ij}&=0\qquad&&\text{if}\quad i>j>k>l\\
s_{ij}\ts s_{kl}-s_{kl}\ts s_{ij}&=0\qquad&&\text{if}\quad i>k>l>j\\
s_{ij}\ts s_{kl}-s_{kl}\ts s_{ij}&=(q-\qin)\tss(s_{kj}s_{il}-s_{ik}s_{jl})
\qquad&&\text{if}\quad i>k>j>l\\
q\ts s_{ij}\ts s_{jl}-s_{jl}\ts s_{ij}&=(q-\qin)\tss s_{il}
\qquad&&\text{if}\quad i>j>l\\
q\ts s_{ij}\ts s_{il}-s_{il}\ts s_{ij}&=(q-\qin)\tss s_{lj}
\qquad&&\text{if}\quad i>l>j\\
q\ts s_{ij}\ts s_{kj}-s_{kj}\ts s_{ij}&=(q-\qin)\tss s_{ki}
\qquad&&\text{if}\quad k>i>j.
\eal
\eeq
In this form the relations were given by Nelson and Regge~\cite{nr:qg}.
An analogue of the Poincar\'e--Birkhoff--Witt theorem for
the algebra $\U'_q(\oa_N)$ was proved in \cite{ik:nd};
see also \cite{m:rtq, mrs:cs}
for other proofs. This theorem implies that at $q=1$ the
algebra $\U'_q(\oa_N)$ specializes to the algebra of
polynomials in $N(N-1)/2$ variables.
More precisely, set $\A=\CC[q,\qin]$ and consider
the $\A$-subalgebra $\U'_{\A}$ of $\U'_q(\oa_N)$
generated by the elements $s_{ij}$.
Then we have an isomorphism
\beql{speci}
\U'_{\A}{\ot}^{}_{\A}\ts\CC\cong \Pc_N,
\eeq
where the action of $\A$ on
$\CC$ is defined via the evaluation $q=1$ and $\Pc_N$ denotes the algebra
of polynomials in the independent variables
$a_{ij}$ with $1\leqslant j<i\leqslant N$. The elements $a_{ij}$
are respective images of the $s_{ij}$ under the isomorphism \eqref{speci}.
Furthermore, the algebra $\Pc_N$ is equipped with the Poisson bracket
$\{\cdot,\cdot\}$
defined by
\beql{defpois}
\{f,h\}=\frac{\wt f\ts\wt h-\wt h\tss\wt f}{1-q}\ts\Big|_{q=1},
\eeq
where $f,h\in\Pc_N$ and $\wt f$ and $\wt h$ are
elements of $\U'_{\A}$ whose images in $\Pc_N$
under the specialization $q=1$ coincide
with $f$ and $h$, respectively.
Indeed, write the element $\wt f\ts\wt h-\wt h\tss\wt f\in\U'_{\A}$
as a linear combination of the ordered monomials in the generators
with coefficients in $\A$.
Since the image of
$\wt f\ts\wt h-\wt h\tss\wt f$ in $\Pc_N$ is zero, all
the coefficients are divisible by $1-q$. Clearly,
the element $\{f,h\}\in \Pc_N$ is independent of the
choice of $\wt f$ and $\wt h$ and of the ordering
of the generators of $\U'_{\A}$.
Obviously, \eqref{defpois} does define
a Poisson bracket on $\Pc_N$.
By definition,
\ben
\{a_{ij},a_{kl}\}=
\frac{s_{ij}\tss s_{kl}-s_{kl}\tss s_{ij}}{1-q}\ts\Big|_{q=1}.
\een
Hence, using the defining relations \eqref{drabssi}, we get
\beql{poissonbr}
\bal
\{a_{ij},a_{kl}\}&=0\qquad&&\text{if}\quad i>j>k>l\\
\{a_{ij},a_{kl}\}&=0\qquad&&\text{if}\quad i>k>l>j\\
\{a_{ij},a_{kl}\}&=2\tss(a_{ik}a_{jl}-a_{kj}a_{il})
\qquad&&\text{if}\quad i>k>j>l\\
\{a_{ij},a_{jl}\}&=a_{ij}a_{jl}-2\tss a_{il}
\qquad&&\text{if}\quad i>j>l\\
\{a_{ij},a_{il}\}&=a_{ij}a_{il}-2\tss a_{lj}
\qquad&&\text{if}\quad i>l>j\\
\{a_{ij},a_{kj}\}&=a_{ij}a_{kj}-2\tss a_{ki}
\qquad&&\text{if}\quad k>i>j.
\eal
\eeq
This coincides with the Poisson brackets of
\cite{b:sg}, \cite{nr:gg}, and \cite{u:ps}, up to a constant factor
if we interpret $a_{ij}$
as the $ji$-th entry of the upper triangular matrix.

We shall also use the presentation of the algebra $\U'_q(\oa_N)$
due to Gavrilik and Klimyk~\cite{gk:qd}. An isomorphism between
the presentations was given by Noumi~\cite{n:ms}, a proof can be
found in Iorgov and Klimyk~\cite{ik:nd}.
Set $s_i=s_{i+1,i}$ for $i=1,\dots,N-1$.
Then the algebra $\U'_q(\oa_N)$ is generated
by the elements $s_1,\dots,s_{N-1}$ subject only to the
relations
\ben
\bal
s^{}_k\tss s_{k+1}^2-(q+\qin)\tss s^{}_{k+1}\tss s^{}_k\tss s^{}_{k+1}
+s_{k+1}^2\tss s^{}_k&=-\qin\tss (q-\qin)^2\tss s^{}_k,\\
s^2_k\tss s_{k+1}^{}-(q+\qin)\tss s^{}_{k}\tss s^{}_{k+1}\tss s^{}_{k}
+s_{k+1}^{}\tss s^2_k&=-\qin\tss (q-\qin)^2\tss s^{}_{k+1},
\eal
\een
for $k=1,\dots,N-2$ (the {\it Serre type relations}), and
\ben
s_k\tss s_l=s_l\tss s_k,\qquad |k-l|>1.
\een

It is easy to see that the subalgebra $\U'_q(\oa_N)\subset\U_q(\gl_N)$
is not preserved by the action of the braid group $B_N$
on $\U_q(\gl_N)$ described above.
Nevertheless, we have the following theorem.

\bth\label{thm:braid}
For $i=1,\dots, N-1$ the assignment
\ben
\bal
\be_i:
    {}&{}s_{i+1}\mapsto \frac{1}{q-\qin}
    \big(\tss q\ts s_{i+1}\ts s_i-s_i\ts s_{i+1}\big)\\
    &s_{i-1}\mapsto \frac{1}{q-\qin}
    \big(s_{i}\ts s_{i-1}-q\ts s_{i-1}\ts s_{i}\big)\\
    &s_{i}\mapsto -s_{i}\\
    &s_{k}\mapsto s_{k}
    \qquad\qquad\text{if}\quad k\ne i-1,i,i+1,
\eal
\een
defines an action of the braid group $B_N$ on $\U'_q(\oa_N)$
by automorphisms.
\eth

\bpf
We verify first that the images
of the generators $s_1,\dots,s_{N-1}$ under $\be_i$ satisfy
the defining relations of $\U'_q(\oa_N)$.
A nontrivial calculation is only required to verify that
the images of the pairs of generators $\be_i(s_k)$ and $\be_i(s_{k+1})$
with $k=i-2,i-1,i,i+1$ satisfy both Serre type relations,
and that the images $\be_i(s_{i-1})$ and $\be_i(s_{i+1})$ commute.
Observe that by \eqref{drabssi},
the image of $s_{i+1}$ can also be written as
\ben
\be_i: s_{i+1}\mapsto s_{i+2,i}.
\een
Hence, for $k=i+1$ we need to verify that
\begin{multline}
s^{}_{i+2,i}\tss s_{i+3,i+2}^2-(q+\qin)\tss s^{}_{i+3,i+2}
\tss s^{}_{i+2,i}\tss s^{}_{i+3,i+2}
+s_{i+3,i+2}^2\tss s^{}_{i+2,i}\\
=-\qin\tss (q-\qin)^2\tss s^{}_{i+2,i}.
\non
\end{multline}
We shall verify the following more general relation in $\U'_q(\oa_N)$,
\beql{genserre}
s^{}_{ij}\tss s_{ki}^2-(q+\qin)\tss s^{}_{ki}
\tss s^{}_{ij}\tss s^{}_{ki}
+s_{ki}^2\tss s^{}_{ij}
=-\qin\tss (q-\qin)^2\tss s^{}_{ij},
\eeq
where $k>i>j$. Indeed,
the left hand side equals
\beql{axilre}
-(q\ts s_{ki}\ts s_{ij}-s_{ij}\ts s_{ki})\ts s_{ki}
+\qin\ts s_{ki}(q\ts s_{ki}\ts s_{ij}-s_{ij}\ts s_{ki}).
\eeq
However, by \eqref{drabssi}
we have
\ben
q\ts s_{ki}\ts s_{ij}-s_{ij}\ts s_{ki}=(q-\qin)\ts s_{kj}
\een
so that \eqref{axilre} becomes
\ben
-\qin\tss (q-\qin)(q\ts s_{kj}\ts s_{ki}-s_{ki}\ts s_{kj})
\een
which equals $-\qin\tss (q-\qin)^2\tss s^{}_{ij}$ by \eqref{drabssi}
thus proving \eqref{genserre}.
The second Serre type relation for the images
$\be_i(s_{i+1})$ and $\be_i(s_{i+2})$ follows from
a more general relation in $\U'_q(\oa_N)$,
\ben
s^2_{ij}\tss s_{ki}^{}-(q+\qin)\tss s^{}_{ij}
\tss s^{}_{ki}\tss s^{}_{ij}
+s_{ki}^{}\tss s^2_{ij}
=-\qin\tss (q-\qin)^2\tss s^{}_{ki},
\een
where $k>i>j$, and which is verified in the same way as \eqref{genserre}.
Next, the Serre type relations for the images
$\be_i(s_{i})$ and $\be_i(s_{i+1})$ follow respectively from
the relations
\ben
s^2_{ij}\tss s_{kj}^{}-(q+\qin)\tss s^{}_{ij}
\tss s^{}_{kj}\tss s^{}_{ij}
+s_{kj}^{}\tss s^2_{ij}
=-\qin\tss (q-\qin)^2\tss s^{}_{ij}
\een
and
\ben
s^2_{ij}\tss s_{kj}^{}-(q+\qin)\tss s^{}_{ij}
\tss s^{}_{kj}\tss s^{}_{ij}
+s_{kj}^{}\tss s^2_{ij}
=-\qin\tss (q-\qin)^2\tss s^{}_{kj},
\een
where $k>i>j$, which both are implied by \eqref{drabssi}.
The Serre type relations for the pairs
$\be_i(s_{i-1}),\ts\be_i(s_{i})$ and $\be_i(s_{i-2}),\ts\be_i(s_{i-1})$
can now be verified by using the involutive automorphism $\om$ of
$\U'_q(\oa_N)$ which is defined on the generators by
\beql{omegade}
s_k\mapsto s_{N-k},\qquad k=1,\dots,N-1.
\eeq
We have
\ben
\bal
\om:{}&{}\be_i(s_{i-2})\mapsto \be_{N-i}(s_{N-i+2}),\\
&\be_i(s_{i-1})\mapsto -\be_{N-i}(s_{N-i+1}),\\
&\be_i(s_{i})\mapsto \be_{N-i}(s_{N-i}),
\eal
\een
and so the desired relations are implied by
the Serre type relations for the pairs of the images
$\be_j(s_{j}),\ts\be_j(s_{j+1})$ and $\be_j(s_{j+1}),\ts\be_j(s_{j+2})$
with $j=N-i$.

Now we verify that the images $\be_i(s_{i-1})$ and $\be_i(s_{i+1})$ commute,
that is,
\beql{commim}
(s_i\ts s_{i-1}-q\ts s_{i-1}\ts s_i)(q\ts s_{i+1}\ts s_i-s_i\ts s_{i+1})
=(q\ts s_{i+1}\ts s_i-s_i\ts s_{i+1})(s_i\ts s_{i-1}-q\ts s_{i-1}\ts s_i).
\eeq
By the Serre type relations we have
\ben
s^2_i\tss s_{i+1}^{}-(q+\qin)\tss s^{}_{i}\tss s^{}_{i+1}\tss s^{}_{i}
+s_{i+1}^{}\tss s^2_i=-\qin\tss (q-\qin)^2\tss s^{}_{i+1}
\een
and
\ben
s^2_i\tss s_{i-1}^{}-(q+\qin)\tss s^{}_{i}\tss s^{}_{i-1}\tss s^{}_{i}
+s_{i-1}^{}\tss s^2_i=-\qin\tss (q-\qin)^2\tss s^{}_{i-1}.
\een
Multiply the first of these relations by $s_{i-1}$
and the second by $s_{i+1}$ from the left. Taking the difference we come to
\ben
s_{i-1}^{}\ts s^2_i\tss s_{i+1}^{}-
(q+\qin)\tss s_{i-1}^{}\ts s^{}_{i}\tss s^{}_{i+1}\tss s^{}_{i}
=s_{i+1}^{}\ts s^2_i\tss s_{i-1}^{}-
(q+\qin)\tss s_{i+1}^{}\ts s^{}_{i}\tss s^{}_{i-1}\tss s^{}_{i}.
\een
Now repeat the same calculation but multiply the Serre type relations
by $s_{i-1}$ and $s_{i+1}$, respectively, from the right.
This gives
\ben
s_{i-1}^{}\ts s^2_i\tss s_{i+1}^{}-
(q+\qin)\tss s_i^{}\ts s_{i-1}^{}\ts s^{}_{i}\tss s^{}_{i+1}
=s_{i+1}^{}\ts s^2_i\tss s_{i-1}^{}-
(q+\qin)\tss s^{}_{i}\tss s^{}_{i+1}\tss s^{}_{i}\ts s_{i-1}^{}.
\een
Hence,
\ben
s_{i-1}^{}\ts s^{}_{i}\tss s^{}_{i+1}\tss s^{}_{i}
-s_{i+1}^{}\ts s^{}_{i}\tss s^{}_{i-1}\tss s^{}_{i}
=s_i^{}\ts s_{i-1}^{}\ts s^{}_{i}\tss s^{}_{i+1}
-s^{}_{i}\tss s^{}_{i+1}\tss s^{}_{i}\ts s_{i-1}^{}
\een
and \eqref{commim} follows.

Thus, each $\be_i$ with $i=1,\dots,N-1$ defines
a homomorphism $\U'_q(\oa_N)\to\U'_q(\oa_N)$. Now observe that
$\be_i$ is invertible with the inverse given by
\ben
\bal
\be^{-1}_i:{}&{}
    s_{i+1}\mapsto \frac{1}{q-\qin}
    \big(s_{i+1}\ts s_i-\tss q\ts s_i\ts s_{i+1}\big)\\
    &s_{i-1}\mapsto \frac{1}{q-\qin}
    \big(\tss q\ts s_{i}\ts s_{i-1}-s_{i-1}\ts s_{i}\big)\\
    &s_{i}\mapsto -s_{i}\\
    &s_{k}\mapsto s_{k}
    \qquad\qquad\text{if}\quad k\ne i-1,i,i+1,
\eal
\een
and so $\be_i$ and $\be^{-1}_i$ are mutually inverse automorphisms
of $\U'_q(\oa_N)$.

Finally, we verify that the automorphisms $\be_i$ satisfy the
braid group relations. It suffices to check that for
each generator $s_k$ we have
\beql{brarelk}
\be_i\tss \be_{i+1}\tss \be_i(s_k)=\be_{i+1}\tss \be_i\tss \be_{i+1}(s_k)
\eeq
for $i=1,\dots,N-2$, and
\beql{commuijk}
\be_i\tss \be_j(s_k)=\be_j\tss \be_i(s_k)
\eeq
for $|i-j|>1$.
Clearly, the only nontrivial cases of \eqref{brarelk} are $k=i-1,i,i+1,i+2$
while \eqref{commuijk} is obvious for all cases except for
$j=i+2$ and $k=i+1$. Take $k=i-1$ in \eqref{brarelk}. We have
$\be_{i+1}(s_{i-1})=s_{i-1}$ while
\ben
\be_{i}:s_{i-1}\mapsto \frac{1}{q-\qin}
\big(s_{i}\ts s_{i-1}-q\ts s_{i-1}\ts s_{i}\big)
=q\ts s_{i+1,i-1}-q\ts s_{i+1,i}\ts s_{i,i-1},
\een
where we have used \eqref{drabssi}. Furthermore, using again \eqref{drabssi},
we find
\begin{multline}
\be_{i+1}\tss \be_i:s_{i-1}\mapsto
q^2\ts s_{i+2,i-1}-q^2\ts s_{i+2,i+1}\ts s_{i+1,i-1}\\
{}-q^2\ts s_{i+2,i}\ts s_{i,i-1}+q^2\ts s_{i+2,i+1}\ts s_{i+1,i}\ts s_{i,i-1}.
\non
\end{multline}
It remains to verify with the use of \eqref{drabssi} that this
element is stable under the action of $\be_i$.
The remaining cases of \eqref{brarelk} and \eqref{commuijk}
are verified with similar and even simpler
calculations.
\epf

\bco\label{cor:braidref}
In terms of the generators $s_{kl}$ of the algebra $\U'_q(\oa_N)$,
for each index $i=1,\dots, N-1$
the action of $\be_i$ is given by
\ben
\bal
\be_i:{}&{}s_{i+1,i}\mapsto -s_{i+1,i}\\
    &s_{ik}\mapsto q\ts s_{i+1,k}-q\ts s_{i+1,i}\tss s_{ik},
    \qquad
    &&s_{i+1,k}\mapsto s_{ik},
    \qquad&&\text{if}\quad k\leqslant i-1\\
    &s_{li}\mapsto \qin\ts s_{l,i+1}-s_{li}\tss s_{i+1,i},
    \qquad
    &&s_{l,i+1}\mapsto s_{li},
    \qquad&&\text{if}\quad l\geqslant i+2\\
    &s_{kl}\mapsto s_{kl}
    &&\text{in all remaining cases.}
\eal
\een
\eco

\bpf
This follows from the defining relations
\eqref{drabssi}. Indeed, the elements $s_{kl}$ can be expressed
in terms of the generators $s_1,\dots,s_{N-1}$ by induction,
using the relations
\beql{indrel}
s^{}_{kl}=\frac{1}{q-\qin}\ts
\big(\tss q\ts s^{}_{kj}\ts s^{}_{jl}-
s^{}_{jl}\ts s^{}_{kj}\big),
\qquad k>j>l.
\eeq
This determines the action of $\be_i$ on the elements $s_{kl}$ and
the formulas are verified by induction.
\epf

\bre\label{rem:dirch}
It is possible to prove that the formulas of Corollary~\ref{cor:braidref}
define an action of the braid group $B_N$ on $\U'_q(\oa_N)$
by automorphisms only using the presentation \eqref{drabssi}.
However, this leads to a slightly longer calculations
as compared with the proof of Theorem~\ref{thm:braid}.
\par
Note also that the universal enveloping algebra $\U(\oa_N)$
can be obtained as a specialization of $\U'_q(\oa_N)$
in the limit $q\to 1$; see \cite{mrs:cs} for a precise formulation.
In this limit the elements
$s_{ij}/(q-\qin)$ with $i>j$
specialize to the generators $F_{ij}$ of $\oa_N$,
where
$F_{ij}=E_{ij}-E_{ji}$. Hence
the action
of $B_N$ on $\U'_q(\oa_N)$ specializes to the action
of the symmetric group $\Sym_N$ on $\U(\oa_N)$ by permutations
of the indices of the $F_{ij}$.
\qed
\ere

The mapping
\eqref{omegade} can also be extended to the entire
algebra $\U'_q(\oa_N)$ as an anti-automorphism. This
is readily verified with the use of the Serre type relations.
We denote this involutive anti-automorphism of $\U'_q(\oa_N)$
by $\omega'$.

\bpr\label{prop:antiom}
The action of $\omega'$ on the generators $s_{kl}$
is given by
\beql{amprs}
\omega':s_{kl}\mapsto s_{N-l+1,\ts N-k+1},\qquad 1\leqslant l<k\leqslant N.
\eeq
Moreover, we have the relations
\beql{ombeom}
\omega'\be^{}_i\ts\omega'=\be^{-1}_{N-i},\qquad i=1,\dots,N-1,
\eeq
where the automorphisms $\be_i$ of $\U'_q(\oa_N)$ are defined
in Theorem~\ref{thm:braid}.
\epr

\bpf
The defining relations \eqref{drabssi} imply that
the mapping \eqref{amprs} defines an anti-automorphism of $\U'_q(\oa_N)$.
Obviously, the images of the generators $s_k$ are found
by \eqref{omegade}. The second part of the proposition
is verified by comparing
the images of the generators $s_k$ under the automorphisms
on both sides of \eqref{ombeom}.
\epf

Observe that the image of the matrix $S$ under $\omega'$ is
given by $\omega':S\mapsto S^{\tss\prime}$, where the prime denotes
the transposition with respect to the second diagonal.

Now consider the involutive automorphism $\omega$ of $\U'_q(\oa_N)$
defined by the mapping \eqref{omegade}.

\bpr\label{prop:autoom}
The image of the matrix $S$ under $\omega$
is given by
\beql{amprsau}
\omega:S\mapsto (1-\qin)\ts I+\qin\ts D(S^{-1})'D^{-1},
\eeq
where $I$ is the identity matrix and
$D=\text{\rm diag}\ts(-q,(-q)^2,\dots,(-q)^{N})$.
In terms of the generators, this can be written as
\ben
\omega:s_{kl}\mapsto (-q)^{k-l-1}
\sum_{N-l+1>r_1>\dots>r_p>N-k+1} (-1)^{p} \ts
s_{N-l+1,r_1}\ts s_{r_1r_2}\dots s_{r_p,N-k+1},\quad k>l,
\een
summed over $p\geqslant 0$ and the indices
$r_1,\dots,r_p$.
\epr

\bpf
The elements $s_{kl}$ can be expressed
in terms of the generators $s_1,\dots,s_{N-1}$
by \eqref{indrel}.
The formula
for $\om(s_{kl})$ is then verified by induction on $k-l$.
The matrix form \eqref{amprsau} is implied by the relation
\beql{sinverse}
(S^{-1})_{kl}=
\sum_{k>r_1>\dots>r_p>l} (-1)^{p+1} \ts
s_{k,r_1}\ts s_{r_1r_2}\dots s_{r_p,l},\qquad k>l,
\eeq
summed over $p\geqslant 0$ and the indices
$r_1,\dots,r_p$.
\epf

For any diagonal matrix $C=\text{diag}\ts(c_1,\dots,c_N)$
the relation \eqref{rsrs}
is preserved by the transformation
$S\mapsto C\tss S\tss C$. Indeed, the entries of $S$ are then
transformed as $s_{ij}\mapsto s_{ij}\ts c_i\ts c_j$ and the claim
is immediate from \eqref{drabs}. This implies that
if $c_i^2=1$ for all $i$ then
the mapping
$\varsigma:S\mapsto C\tss S\tss C$ defines an automorphism of $\U'_q(\oa_N)$.
Therefore, Propositions~\ref{prop:antiom} and \ref{prop:autoom}
imply the following corollary.

\bco\label{cor:sinv}
The mapping
\beql{amprsauant}
\rho:S\mapsto (1-\qin)\ts I+\qin\ts H\tss S^{-1}\tss H^{-1},
\eeq
where
$H=\text{\rm diag}\ts(q,q^{2},\dots,q^N)$,
defines an involutive anti-automorphism
of $\U'_q(\oa_N)$.
\eco

\bpf
We obviously have $\rho=\varsigma\circ\omega'\circ\omega$
for an appropriate automorphism $\varsigma$. Hence
$\rho$ is an anti-automorphism. We have
\ben
\rho:s_k\mapsto -s_k,\qquad k=1,\dots,N-1,
\een
and so $\rho$ is involutive.
\epf

We can now recover the braid group action on the algebra
$\Pc_N$; see Dubrovin~\cite{d:gt}, Bondal~\cite{b:sg}.

\bco\label{cor:braidpol}
The braid group $B_N$ acts on the algebra $\Pc_N$ by
\ben
\bal
\be_i:{}&{}a_{i+1,i}\mapsto -a_{i+1,i}\\
    &a_{ik}\mapsto a_{i+1,k}-a_{i+1,i}\tss a_{ik},
    \qquad
    &&a_{i+1,k}\mapsto a_{ik},
    \qquad&&\text{if}\quad k\leqslant i-1\\
    &a_{li}\mapsto a_{l,i+1}-a_{li}\tss a_{i+1,i},
    \qquad
    &&a_{l,i+1}\mapsto a_{li},
    \qquad&&\text{if}\quad l\geqslant i+2\\
    &a_{kl}\mapsto a_{kl}
    &&\text{in all remaining cases},
\eal
\een
where $i=1,\dots,N-1$. Moreover, the Poisson bracket on $\Pc_N$
in invariant under this action.
\eco

\bpf
This is immediate from Corollary~\ref{cor:braidref}.
\epf

We combine the variables $a_{ij}$ into the lower triangular matrix
$A=[a_{ij}]$ where we set $a_{ii}=1$ for all $i$ and $a_{ij}=0$
for $i<j$.

\bco\label{cor:poisom}
The mappings
\beql{amprsaua}
\varrho:A\mapsto A^{-1}\Fand \varpi:A\mapsto A'
\eeq
define anti-automorphisms of
the Poisson bracket on $\Pc_N$.
Explicitly, the image of $a_{kl}$ under $\varrho$ is given by
\ben
\varrho:a_{kl}\mapsto
\sum_{k>r_1>\dots>r_p>l} (-1)^{p+1} \ts
a_{kr_1}\ts a_{r_1r_2}\dots a_{r_p,l},\quad k>l,
\een
summed over $p\geqslant 0$ and the indices
$r_1,\dots,r_p$.
\eco

\bpf
This follows from Proposition~\ref{prop:antiom}
and Corollary~\ref{cor:sinv} by taking $q=1$.
\epf

\section{Casimir elements of the Poisson algebra $\Pc_N$}\label{sec:cep}
\setcounter{equation}{0}

Using the relationship between the twisted
quantized enveloping algebra $\U'_q(\oa_N)$ and the
Poisson algebra $\Pc_N$, we can get families of invariants
of $\Pc_N$ by taking the classical limit $q\to 1$
in the constructions of \cite{mrs:cs}, \cite{gi:ce} and \cite{nuw:dp}.
First, we recall the construction of
Casimir elements for
the algebra $\U'_q(\oa_N)$ given in \cite{mrs:cs}.
Consider the $q$-permutation operator
$P^{\tss q}\in\End(\CC^N\ot\CC^N)$
defined by
\beql{qperm}
P^{\tss q}=\sum_{i}E_{ii}\ot E_{ii}+ q\tss\sum_{i> j}E_{ij}\ot
E_{ji}+ \qin\sum_{i< j}E_{ij}\ot E_{ji}.
\end{equation}
Introduce
the multiple tensor product
$\U'_q(\oa_N)\ot (\End\CC^N)^{\ot\tss r}$.
The action of the symmetric group $\Sym_r$ on the space $(\CC^N)^{\ot\tss r}$
can be defined by setting $\si_i\mapsto P^{\tss q}_{\si_i}:=
P^{\tss q}_{i,i+1}$ for $i=1,\dots,r-1$,
where $\si_i$ denotes the transposition $(i,i+1)$.
If $\si=\si_{i_1}\cdots \si_{i_l}$ is a reduced decomposition
of an element $\si\in \Sym_r$ we set
$P^{\tss q}_{\si}=P^{\tss q}_{\si_{i_1}}\cdots P^{\tss q}_{\si_{i_l}}$.
We denote by $A^q_r$ the $q$-antisymmetrizer
\beql{antisym}
A^q_r=\sum_{\sigma\in\Sym_r}\sgn\tss\si\cdot P^{\tss q}_{\si}.
\eeq
Now take $r=N$. We have the relation
\beql{sdetmatr}
\bal
A^q_N\ts S_1(u)\ts R_{12}^{\tss t}\cdots R_{1N}^{\tss t}\ts S_2(u\tss q^{-2})
\ts &R_{23}^{\tss t}\cdots R_{2N}^{\tss t}\ts
S_3(u\tss q^{-4})\\
&{}\times \cdots R_{N-1,N}^{\tss t}\ts
S_N(u\tss q^{-2N+2})\\[0.5em]
{}=S_N(u\tss q^{-2N+2})\ts R_{N-1,N}^{\tss t}
\cdots S_3(u\tss q^{-4})&\ts R_{2N}^{\tss t}\cdots R_{23}^{\tss t}\ts
S_2(u\tss q^{-2})\\
&{}\times\ts R_{1N}^{\tss t} \cdots R_{12}^{\tss t}
\ts  S_1(u)\ts A^q_N,
\eal
\eeq
where the following notation was used. The matrix $S(u)$
is defined by
\ben
S(u)=S+\qin\ts u^{-1}\ts {\overline S},
\een
where $u$ is a formal
variable and
${\overline S}$ is the upper triangular matrix with ones on the diagonal
whose $ij$-th entry is $\bar s_{ij}=q\ts s_{ji}$
for $i<j$. Furthermore,
\ben
R_{ij}^{\tss t}=R_{ij}^{\tss t}(u^{-1}\tss q^{\tss 2i-2},u\tss q^{-2j+2})
\een
with
\beql{rtuv}
\bal
R^{\ts t}(u,v)={}&(u-v)\sum_{i\ne j}E_{ii}\ot E_{jj}+(\qin u-q\tss v)
\sum_{i}E_{ii}\ot E_{ii} \\
{}+ {}&(\qin-q)\tss u\tss\sum_{i> j}E_{ji}\ot
E_{ji}+ (\qin-q)\tss v\tss\sum_{i< j}E_{ji}\ot E_{ji}.
\eal
\eeq
The subscripts in \eqref{sdetmatr} indicate the copies
of $\End\CC^N$ in
$\U'_q(\oa_N)\ot (\End\CC^N)^{\ot\tss N}$ which are labelled by
$1,\dots,N$; cf. \eqref{defrel}. The element \eqref{sdetmatr}
equals $A^q_N\ts\sdet S(u)$, where $\sdet S(u)$ is a rational
function
in $u$ (the {\it Sklyanin determiant\/})
valued in the center of $\U'_q(\oa_N)$; see
\cite[Theorem~3.8 and Corollary~4.3]{mrs:cs}.

Recall that the Poisson algebra $\Pc_N$ is the algebra of
polynomials in the variables $a_{ij}$ with $i>j$
which are combined into the matrix
$A=[a_{ij}]$ with $a_{ii}=1$ for all $i$ and $a_{ij}=0$
for $i<j$.
The following theorem was proved in different ways
by Nelson and Regge~\cite{nr:ig}
and Bondal~\cite{b:sg}.

\bth\label{thm:casimirs}
The coefficients of the polynomial
\ben
\det(A+\la A^t)=f_0+f_1\tss\la+\dots+f_N\tss\la^N
\een
are Casimir elements of the Poisson algebra $\Pc_N$.
\eth

\bpf
We use the centrality of the Sklyanin determinant $\sdet S(u)$
in $\U'_q(\oa_N)$.
Note that at $q=1$ the $q$-antisymmetrizer $A^q_N$ becomes
the antisymmetrizer in $(\CC^N)^{\ot\tss N}$, the element
$R^{\tss t}(u-v)$ becomes $u-v$ times the identity. Since the images of
the elements $s_{ij}$ in $\Pc_N$ coincide with $a_{ij}$,
the image of the matrix $S(u)$ is $A+u^{-1} A^t$.
Hence, at $q=1$ the Sklyanin determinant $\sdet S(u)$
becomes $\ga(u)\tss\det(A+u^{-1} A^t)$, where
\beql{gau}
\ga(u)=(u^{-1}-u)^{N(N-1)/2}.
\eeq
Therefore, replacing $u$ with $\la^{-1}$
we thus prove that
all coefficients
of $\det(A+\la A^t)$ are Casimir elements for the Poisson
bracket on $\Pc_N$.
\epf

Note that, as was proved in \cite{b:sg} and \cite{nr:ig},
the polynomial $\det(A+\la A^t)$ is invariant under the action
of the braid group $B_N$.

Now we recall the construction of Casimir elements
given in \cite{gi:ce}. For all $i>j$
define the elements $s^+_{ij}$
of $\U'_q(\oa_N)$ by induction from the formulas
\ben
s^{+}_{ij}=\frac{1}{q-\qin}\ts
\big(s^{+}_{i,\tss j+1}\ts s^{}_{j+1,\tss j}-
q\ts s^{}_{j+1,\tss j}\ts s^{+}_{i,\tss j+1}\big),
\qquad i>j+1,
\een
and $s^{+}_{j+1,\tss j}=s^{}_{j+1,\tss j}$ for $j=1,\dots,N-1$.
A straightforward calculation shows that
these elements can be equivalently defined by
\ben
s^{+}_{ij}=-q^{i-j-1}\ts (S^{-1})_{ij},\qquad i>j,
\een
where the entries of the inverse matrix are found from
\eqref{sinverse}.
Let $k$ be a positive integer such that $2k\leqslant N$.
For any subset $I=\{i_1<i_2<\dots<i_{2k}\}$
of $\{1,\dots,N\}$ introduce the elements $\Phi^{}_I$ and $\Phi_I^+$
of $\U'_q(\oa_N)$ by
\ben
\Phi^{}_I=\sum_{\si\in\Sym_{2k}}(-q)^{-\ell(\si)}
s^{}_{i_{\si(2)}\ts i_{\si(1)}}\dots s^{}_{i_{\si(2k)}\ts i_{\si(2k-1)}}
\een
and
\ben
\Phi^{+}_I=\sum_{\si\in\Sym_{2k}}(-q)^{\ell(\si)}
s^+_{i_{\si(2)}\ts i_{\si(1)}}\dots s^+_{i_{\si(2k)}\ts i_{\si(2k-1)}},
\een
where $\ell(\si)$ is the length of the permutation $\si$, and the sums
are taken over those permutations $\si\in\Sym_{2k}$ which satisfy
the conditions
\ben
i_{\si(2)}> i_{\si(1)},\quad\dots,\quad i_{\si(2k)}> i_{\si(2k-1)}
\Fand i_{\si(2)}<i_{\si(4)}<\dots<i_{\si(2k)}.
\een
Then according to \cite{gi:ce}, for each $k$ the element
\ben
\phi_k=\sum_{I,\ts |I|=2k} q^{i_1+i_2+\dots+i_{2k}}_{}
\ts\Phi^{+}_I\Phi^{}_I
\een
belongs to the center of $\U'_q(\oa_N)$. Moreover,
in the case $N=2n$ both elements $\Phi^{}_{I_0}$ and $\Phi^{+}_{I_0}$
with $I_0=\{1,\dots,2n\}$ are also central.

\bre\label{rem:ortlim}
Our notation is related to \cite{gi:ce} by
\ben
s^{}_{ij}=-q^{-1/2}(q-\qin)\ts I^-_{ij},\qquad
s^+_{ij}=-q^{-1/2}(q-\qin)\ts I^+_{ij},\qquad i>j.
\een
Note also that
the elements $\phi_k$ are $q$-analogues of the Casimir elements
for the orthogonal Lie algebra $\oa_N$ constructed in \cite{mn:ci};
see also \cite{iu:ce}.
\qed
\ere

Now return to the Poisson algebra $\Pc_N$.
Recall that the Pfaffian of a $2k\times 2k$ skew symmetric matrix $H$
is given by
\ben
\Pf\ts H=\frac{1}{2^k\ts k!}\sum_{\si\in\Sym_{2k}}\sgn\si\cdot
H_{\si(1),\si(2)}\dots H_{\si(2k-1),\si(2k)}.
\een
Given a lower triangular
$N\times N$ matrix $B$ and a $2k$-element subset $I$
of $\{1,\dots,N\}$ as above, we denote by
$\Pf^{}_I(B)$ the Pfaffian of the $2k\times 2k$
submatrix $(B^t-B)_I$ of $B^t-B$ whose
rows and columns are determined by the elements
of $I$.

\bth\label{thm:pfaff}
For each positive integer $k$ such that $2k\leqslant N$
the element
\beql{defck}
c_k=(-1)^k\ts\sum_{I,\ts |I|=2k}
\ts\Pf^{}_I(A)\ts\Pf^{}_I(A^{-1})
\eeq
is a Casimir element of $\Pc_N$. Moreover,
in the case $N=2n$ both $\Pf^{}_{I_0}(A)$ and $\Pf^{}_{I_0}(A^{-1})$
with $I_0=\{1,\dots,2n\}$ are also Casimir elements.
\eth

\bpf Observe that
in the limit $q\to 1$ the elements $\Phi^{}_I$ and $\Phi^{+}_I$
specialize respectively to the Pfaffians
\ben
\Phi^{}_I\to \Pf^{}_I(A),\qquad \Phi^{+}_I\to (-1)^k\ts\Pf^{}_I(A^{-1}).
\een
Hence, the central element $\phi_k$ specializes to $c_k$.
\epf

\bex\label{ex:markov}
As the matrix elements of the inverse matrix $A^{-1}$ are found
by the formula of Corollary~\ref{cor:poisom}, we have the following explicit
formula for $c_1$,
\ben
c_1=\sum_{i>r_1>\dots>r_p>j}(-1)^{p} \ts
a_{ij}\ts a_{ir_1}\ts a_{r_1r_2}\dots a_{r_pj}.
\een
For $N=3$ it gives the Markov polynomial.
\qed
\eex

\bco\label{cor:pfgen}
The algebra of Casimir elements of $\Pc_N$
is generated by $c_1,\dots,c_n$ for $N=2n+1$,
and by $c_1,\dots,c_{n-1},\Pf^{}_{I_0}(A)$ if $N=2n$.
In both cases, the families of generators are algebraically
independent. Moreover, $\Pf^{}_{I_0}(A^{-1})=(-1)^n\ts\Pf^{}_{I_0}(A)$.
\eco

\bpf
Since
\ben
\det(A+\la A^t)=\la^N\ts \det(A+\la^{-1} A^t),
\een
we have the relations
$f_{N-i}=f_i$.
Moreover, $f_0=f_N=1$ since $\det A=1$.
It was proved in \cite{b:sg} that if $N=2n+1$ is odd then
the coefficients $f_1,\dots,f_n$ are algebraically independent
generators of the algebra of Casimir elements of $\Pc_N$.
If $N=2n$ is even then
\beql{depfsq}
\det(A-A^t)=\Pf^{}_{I_0}(A)^2.
\eeq
In this case, a family of algebraically independent
generators of the algebra of Casimir elements of $\Pc_N$
is obtained by replacing any one of the elements $f_1,\dots,f_n$
with $\Pf^{}_{I_0}(A)$.
The claim will be implied by
the following identity
\beql{idens}
\det(A+\la A^t)=\sum_{k=0}^n(-\la)^k(1+\la)^{N-2k}\ts c_k.
\eeq
Indeed, by the identity, the elements $f_1,\dots,f_n$ can be expressed
as linear combinations of $c_1,\dots,c_n$. In order to verity \eqref{idens},
we use the observation of \cite{b:sg} that
the Casimir elements of $\Pc_N$ are determined by their
restrictions on a certain subspace $\Hc$ of matrices.
If $N=2n$ then $\Hc$ consists of the matrices of the form
\beql{matfo}
\begin{pmatrix} I&O\\
                D&I
\end{pmatrix},
\eeq
where $I$ and $O$ are the identity and zero $n\times n$
matrices, respectively, while $D=\text{diag}(d_1,\dots,d_n)$
is an arbitrary diagonal matrix. If $N=2n+1$ then
$\Hc$ consists of the matrices obtained from \eqref{matfo}
by inserting an extra row and column in the middle of the matrix
whose only nonzero entry is $1$ at their intersection.
So, by Theorems~\ref{thm:casimirs} and \ref{thm:pfaff},
we only need to verify \eqref{idens} for the matrices $A\in\Hc$.
However, in this case the element $c_k$ coincides
with the elementary symmetric polynomial
\ben
c_k=\sum_{r_1<\dots<r_k} d_{r_1}^2\dots d_{r_k}^2,
\een
while
\ben
\det(A+\la A^t)=\prod_{i=1}^n \big((1+\la)^2-\la\ts d_i^2\big)
\een
if $N=2n$, and
\ben
\det(A+\la A^t)=(1+\la)\prod_{i=1}^n \big((1+\la)^2-\la\ts d_i^2\big)
\een
if $N=2n+1$.
This gives \eqref{idens}. To verify the last statement of the corollary,
put $\la=-1$ into \eqref{idens} with $N=2n$. Together with \eqref{depfsq}
this gives $c_n=\Pf^{}_{I_0}(A)^2$, so that the statement follows from
\eqref{defck} with $k=n$.
\epf

Finally, we consider the invariants of the Poisson bracket on $\Pc_N$
which can obtained from the construction of the Casimir elements
of $\U'_q(\oa_N)$ given in \cite{nuw:dp}.

\bth\label{thm:gelfinv}
The elements
\ben
\text{\rm tr}\ts (A^{-1}A^t)^k,\qquad k=1,2,\dots,
\een
are Casimir elements of $\Pc_N$.
\eth

\bpf
This follows by taking the classical limit of the Casimir elements
of \cite{nuw:dp}. Alternatively, this is also implied by
Theorem~\ref{thm:casimirs} and the Liouville formula
\ben
\sum_{k=1}^{\infty}(-1)^{k-1}\la^{k-1} \text{\rm tr}\ts H^k=
\frac{d}{d\la}\ln\det(1+\la H)
\een
which holds for any square matrix $H$. We apply it to the matrix
$H=A^{-1}A^t$ and observe that $\det(A+\la A^t)=\det(1+\la H)$
since $\det A=1$.
\epf

\section{A new Poisson algebra}\label{sec:npa}
\setcounter{equation}{0}

Here we use the symplectic version of the twisted quantized
enveloping algebra introduced
by Noumi~\cite{n:ms} to define a new Poisson algebra and calculate
its Casimir elements.

The {\it twisted
quantized enveloping algebra\/} $\U'_q(\spa_{2n})$
is an associative algebra generated by elements $s_{ij}$,
$i,j\in\{1,\dots,2n\}$ and $s_{i,i+1}^{-1}$, $i=1,3,\dots,2n-1$.
The generators $s_{ij}$ are zero
for $j=i+1$ with even $i$, and
for $j\geqslant i+2$ and all $i$. We combine the $s_{ij}$ into
a matrix $S$ as in \eqref{matrt},
\beql{matrts}
S=\sum_{i,j}s_{ij}\ot E_{ij},
\eeq
so that $S$ has a block-triangular form with
$n$ diagonal $2\times 2$-blocks,
\ben
S=\begin{pmatrix} s_{11}&s_{12}&0&0&\cdots&0&0\\
                     s_{21}&s_{22}&0&0&\cdots&0&0\\
                     s_{31}&s_{32}&s_{33}&s_{34}&\cdots&0&0\\
                     s_{41}&s_{42}&s_{43}&s_{44}&\cdots&0&0\\
                         \vdots&\vdots&\vdots&\vdots&\ddots&\vdots&\vdots\\
                     s_{2n-1,1}&s_{2n-1,2}
                         &s_{2n-1,3}&s_{2n-1,4}&
                                \cdots&s_{2n-1,2n-1}&s_{2n-1,2n}\\
                     s_{2n,1}&s_{2n,2}&s_{2n,3}&s_{2n,4}&
                                               \cdots&s_{2n,2n-1}&s_{2n,2n}
\end{pmatrix}.
\een
The defining relations of $\U'_q(\spa_{2n})$ have the form
of a reflection equation \eqref{rsrs}
together with
\beql{invrel}
s^{}_{i,i+1}\ts s_{i,i+1}^{-1}=s_{i,i+1}^{-1}\ts s^{}_{i,i+1}=1
\eeq
and
\beql{qdetrel}
s^{}_{i+1,i+1}\ts s^{}_{ii}-q^2\ts s^{}_{i+1,i}s^{}_{i,i+1}=q^3
\eeq
for $i=1,3,\dots,2n-1$. More explicitly, the relations
\eqref{rsrs} have exactly the same form \eqref{drabs}
as in the orthogonal case.

Recall the quantized enveloping algebra $\U_q(\gl_{2n})$
defined in Section~\ref{sec:bga}.
Introduce the block-diagonal $2n\times 2n$ matrix $G$ by
\ben
G=\left(\begin{matrix}
                     0&q&\cdots&0&0\\
                     -1&0&\cdots&0&0\\
                     \vdots&\vdots&\ddots&\vdots&\vdots\\
                     0&0&\cdots&0&q\\
                     0&0&\cdots&-1&0
\end{matrix}\right).
\een
We can regard $\U'_q(\spa_{2n})$ as a subalgebra
of $\U_q(\gl_{2n})$ by setting
$S= T\ts G\ts \overline T^{\ts t}$,
or in terms of generators,
\beql{sijtij}
s_{ij}= q\ts \sum_{k=1}^n
t_{i,2k-1}\ts\bar t_{j,2k}
-\sum_{k=1}^n t_{i,2k}\ts\bar t_{j,2k-1};
\eeq
see \cite{n:ms} and \cite{mrs:cs}
for the proofs.

Define the {\it extended twisted
quantized enveloping algebra\/} $\hat\U'_q(\spa_{2n})$
as follows.
This is an associative algebra
generated by elements $s_{ij}$,
$i,j\in\{1,\dots,2n\}$
where
$s_{ij}=0$
for $j=i+1$ with even $i$, and
for $j\geqslant i+2$ and all $i$.
The defining relations are given by \eqref{rsrs}
or, equivalently, by \eqref{drabs}.
We use the same symbols as for the generators of $\U'_q(\spa_{2n})$;
a confusion should be avoided as we indicate
which algebra is considered
at any moment.
This definition essentially coincides
with the original one due to Noumi~\cite{n:ms}.
Note that, in comparison with $\U'_q(\spa_{2n})$,
we neither require the elements
$s^{}_{i,i+1}$ with odd $i$ to be invertible, nor we impose
the relations \eqref{qdetrel}.

An analogue of the Poincar\'e--Birkhoff--Witt theorem
for the algebra $\hat\U'_q(\spa_{2n})$ follows from
\cite[Corollary~3.4]{m:rtq}.
As with the algebra $\U'_q(\oa_N)$, this theorem implies
that at $q=1$ the extended twisted
quantized enveloping algebra $\hat\U'_q(\spa_{2n})$
specializes to the algebra $\hat\Pc_{2n}$ of
polynomials in $2n^2+2n$ variables. We denote the variables by $a_{ij}$
with the same restrictions on the indices $i,j$ as for the elements $s_{ij}$,
so that $s_{ij}$ specializes to $a_{ij}$. We shall combine
the variables $a_{ij}$ into a matrix $A$ which
has a block-triangular form with
$n$ diagonal $2\times 2$-blocks,
\ben
A=\begin{pmatrix} a_{11}&a_{12}&0&0&\cdots&0&0\\
                     a_{21}&a_{22}&0&0&\cdots&0&0\\
                     a_{31}&a_{32}&a_{33}&a_{34}&\cdots&0&0\\
                     a_{41}&a_{42}&a_{43}&a_{44}&\cdots&0&0\\
                         \vdots&\vdots&\vdots&\vdots&\ddots&\vdots&\vdots\\
                     a_{2n-1,1}&a_{2n-1,2}
                         &a_{2n-1,3}&a_{2n-1,4}&
                                \cdots&a_{2n-1,2n-1}&a_{2n-1,2n}\\
                     a_{2n,1}&a_{2n,2}&a_{2n,3}&a_{2n,4}&
                                               \cdots&a_{2n,2n-1}&a_{2n,2n}
\end{pmatrix}.
\een

\bth\label{thm:poissymp}
The algebra $\hat\Pc_{2n}$ possesses the Poisson bracket
defined by
\ben
\bal
\{a_{ij}, a_{kl}\}&=\big(\tss\de_{ik}+\de_{jk}-\de_{il}-\de_{jl}\big)
\ts a_{ij}\ts a_{kl}\\
{}&-2\tss \big(\tss\de_{l<j} -\de_{i<k}\big)
\ts a_{kj}\ts a_{il}-2\tss \de_{l<i}\ts a_{ki}\ts a_{lj}
+2\tss \de_{j<k}\ts a_{ik}\ts a_{jl}.
\eal
\een
\eth

\bpf
We define the Poisson bracket on $\hat\Pc_{2n}$
by the same rule \eqref{defpois}
as in the orthogonal case. The explicit formulas
for the values $\{a_{ij},a_{kl}\}$ follow from \eqref{drabs}.
\epf

\bre\label{rem:pb}
Both in the orthogonal and symplectic case, the Poisson
brackets of $\Pc=\Pc_N$ or $\Pc=\hat\Pc_{2n}$ can be written
in a uniform way in a matrix form. Introducing
the elements of $\Pc\ot\End\CC^N\ot\End\CC^N$ by
\ben
A_1=\sum_{i,j}a_{ij}\ot E_{ij}\ot 1,\qquad
A_2=\sum_{i,j}a_{ij}\ot 1\ot E_{ij},
\een
we have
\ben
\{A_1,A_2\}=[r,A_1A_2]+A_1\,r^{t}A_2-A_2\,r^{t}A_1,
\een
where
\ben
r=\sum_i    E_{ii}\otimes E_{ii} +2\,\sum_{i<j} E_{ij}\otimes E_{ji},
\qquad
r^t=\sum_i    E_{ii}\otimes E_{ii} +2\,\sum_{i<j} E_{ji}\otimes E_{ji}.
\een
This follows from \eqref{rsrs} and the observation that
\ben
r=\frac{R-I\ot I}{q-1}\Big|^{}_{q=1}.
\een
\ere

\bth\label{thm:casimsym}
The elements
\beql{smdets}
a^{}_{i+1,i+1}\ts a^{}_{ii}-a^{}_{i+1,i}\ts a^{}_{i,i+1},\qquad i=1,3,\dots,2n-1,
\eeq
and the coefficients of the polynomial
\ben
\det(A+\la A^t)=f_0+f_1\tss\la+\dots+f_{2n}\tss\la^{2n}
\een
are Casimir elements of the Poisson algebra $\hat\Pc_{2n}$.
\eth

\bpf
For any $i=1,3,\dots,2n-1$ the element
\ben
s^{}_{i+1,i+1}\ts s^{}_{ii}-q^2\ts s^{}_{i+1,i}\ts s^{}_{i,i+1}
\een
belongs to the center of the algebra
$\hat\U'_q(\spa_{2n})$; see \cite[Section~2.2]{mrs:cs}. This implies
the claim for the elements \eqref{smdets}.

We proceed as in the proof of Theorem~\ref{thm:casimirs}. The relation
\eqref{sdetmatr} holds in the same form with the matrix $S(u)$
now given by
\ben
S(u)=S+q\ts u^{-1}\ts {\overline S},
\een
where the matrix elements $\bar s_{ij}$ of
${\overline S}$ are defined as follows.
For any $i=1,3,\dots, 2n-1$ we have
\ben
\bal
\bar s_{ii}&=-q^{-2}\ts s_{ii},\qquad \bar s_{i+1,i+1}=-q^{-2}\ts s_{i+1,i+1},\\
\bar s_{i+1,i}&=-q^{-1}\ts s_{i,i+1},\qquad
\bar s_{i,i+1}=-q^{-1}\ts s_{i+1,i}+(1-q^{-2})\ts s_{i,i+1},
\eal
\een
while
\ben
\bar s_{kl}=-q^{-1}\ts s_{lk}
\een
for $k<l$ except for the pairs $k=i$, $l=i+1$, with odd $i$, and
the remaining entries of ${\overline S}$ are equal to zero.
The element \eqref{sdetmatr}
equals $A^q_N\ts\sdet S(u)$, where $\sdet S(u)$
is the {\it Sklyanin determiant\/} of the matrix $S(u)$.
This is a rational function in $u$
valued in the (extended) twisted quantized enveloping algebra.
When the values are considered in the algebra $\U'_q(\spa_{2n})$,
they are contained in
the center of $\U'_q(\spa_{2n})$, as proved in
\cite[Theorem~3.15 and Corollary~4.3]{mrs:cs}.
The same property holds for
the algebra
$\hat\U'_q(\spa_{2n})$, that is,
when the values of the function
$\sdet S(u)$ are regarded as elements of
the extended algebra $\hat\U'_q(\spa_{2n})$, they belong to
the center of $\hat\U'_q(\spa_{2n})$
(see the proof in the Appendix).

At $q=1$ the matrix $S(u)$ becomes $A-u^{-1} A^t$.
Hence, the Sklyanin determinant $\sdet S(u)$
becomes $\ga(u)\det(A-u^{-1} A^t)$, where $\ga(u)$
is defined in \eqref{gau} with $N=2n$.
Therefore, replacing $u$ with $-\la^{-1}$
we thus prove that
all coefficients
of $\det(A+\la A^t)$ are Casimir elements for the Poisson
bracket on $\hat\Pc_{2n}$.
\epf

As in the orthogonal case, we have
$f_{2n-i}=f_i$ for all $i=0,1,\dots,2n$.
Note also that $f_0=f_{2n}=\det A$ and so
we have the following relation between the Casimir elements
\ben
f_0=\prod_{k=1}^n
\big(a^{}_{2k,2k}\ts a^{}_{2k-1,2k-1}-a^{}_{2k,2k-1}\ts a^{}_{2k-1,2k}\big).
\een

\bcj\label{conj:casim}
The algebra of Casimir elements of $\hat\Pc_{2n}$
is generated by the family of elements
provided by Theorem~\ref{thm:casimsym}
and the Pfaffian $\Pf(A-A^t)$.
\ecj

In the rest of this section we work with the twisted
quantized enveloping algebra $\U'_q(\spa_{2n})$.
Recall the action of the braid group $B_{2n}$
on the quantized enveloping algebra $\U_q(\gl_{2n})$; see Section~\ref{sec:bga}.

\bpr\label{prop:bodd}
The subalgebra $\U'_q(\spa_{2n})\subset \U_q(\gl_{2n})$ is
stable under the action of the elements
$\be_1,\be_3,\dots,\be_{2n-1}$ of $B_{2n}$.
\epr

\bpf
Observe that the algebra $\U'_q(\spa_{2n})$ is generated by the
elements
\beql{diabl}
s^{}_{ii},\quad s^{}_{i+1,i+1},\quad s^{}_{i,i+1},\quad s_{i,i+1}^{-1}
\qquad \text{for}\quad i=1,3,\dots,2n-1
\eeq
and
\beql{sunddi}
s_{i+3,i+1}\qquad
\text{for}\quad i=1,3,\dots,2n-3.
\eeq
Indeed, $s_{i+1,i}$ for odd $i$ can be expressed in terms of the
elements \eqref{diabl} from \eqref{qdetrel}.
Furthermore, the remaining generators can be expressed
in terms of the elements
\eqref{diabl} and
\beql{undial}
s_{i+2,i},\quad s_{i+2,i+1},\quad s_{i+3,i}\quad s_{i+3,i+1}\qquad
\text{for}\quad i=1,3,\dots,2n-3
\eeq
by induction
from the relations
\ben
(q-\qin)\ts s_{kl}=s_{i,i+1}^{-1}\ts
(s_{k,i+1}\ts s_{il}-s_{il}\ts s_{k,i+1}),\qquad k>i+1,\quad i>l,\quad i\
\text{odd},
\een
which are implied by the defining relations \eqref{drabs}.
However, for each $i$ as in \eqref{undial} we have
\ben
\bal
(q-\qin)\ts s_{i+3,i}&=s_{i,i+1}^{-1}\ts
(s_{i+3,i+1}\ts s_{ii}-s_{ii}\ts s_{i+3,i+1}),\\
(q-\qin)\ts s_{i+2,i+1}&=s_{i+3,i+2}^{-1}\ts
(s_{i+3,i+1}\ts s_{i+2,i+2}-s_{i+2,i+2}\ts s_{i+3,i+1}),\\
(q-\qin)\ts s_{i+2,i}&=s_{i+3,i+2}^{-1}\ts
(s_{i+3,i}\ts s_{i+2,i+2}-s_{i+2,i+2}\ts s_{i+3,i}).
\eal
\een

Hence, it suffices to verify that the images of the elements
\eqref{diabl} and \eqref{sunddi} under
the action of $\be_1,\be_3,\dots,\be_{2n-1}$
are contained in $\U'_q(\spa_{2n})$. These images
can be explicitly calculated from \eqref{sijtij}. For any odd $j$
the elements \eqref{diabl} with $i\ne j$ are fixed by the action
of $\be_j$, while
\ben
\be_j: s^{}_{jj}\mapsto s_{j,j+1}^{-2}\ts s^{}_{j+1,j+1},
\qquad s^{}_{j+1,j+1}\mapsto q^{-2}\ts s^{}_{jj},\qquad
s^{}_{j,j+1}\mapsto q^{2} s_{j,j+1}^{-1}.
\een
Moreover, the elements \eqref{sunddi} with $i\ne j-2,j$ are fixed by the action
of $\be_j$, while
\ben
\be_j: s^{}_{j+1,j-1}\mapsto \qin\ts s^{}_{j,j-1},\qquad
s^{}_{j+3,j+1}\mapsto \qin\ts s^{}_{j+3,j}.
\een
All these relations are verified by direct calculation
with the use of the defining relations of $\U_q(\gl_{2n})$.
\epf

In particular, the restrictions of the action of
$\be_1,\be_3,\dots,\be_{2n-1}$
to the subalgebra $\U'_q(\spa_{2n})$ yield automorphisms of the latter.

Now observe that
the elements $\ga_1,\ga_3,\dots,\ga_{2n-3}$ of $B_{2n}$
given by
\ben
\ga_{2k-1}=\be_{2k}\be_{2k-1}\be_{2k+1}\be_{2k},\qquad k=1,\dots,n-1
\een
generate a subgroup of $B_{2n}$ isomorphic to $B_n$. The braid relations
for the $\ga_{2k-1}$ are easily verified with the use of their
geometric interpretation.
Indeed, if we regard $\be_j$ as the braid

\vspace{-2em}

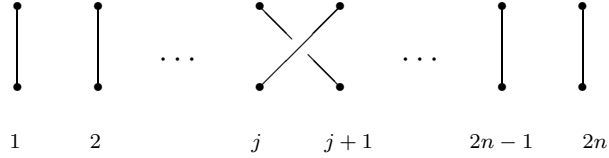
\setlength{\unitlength}{0.13em}
\begin{center}
\begin{picture}(200,60)

\thinlines

\put(10,20){\circle*{2}}
\put(30,20){\circle*{2}}
\put(70,20){\circle*{2}}
\put(90,20){\circle*{2}}
\put(130,20){\circle*{2}}
\put(150,20){\circle*{2}}

\put(10,40){\circle*{2}}
\put(30,40){\circle*{2}}
\put(70,40){\circle*{2}}
\put(90,40){\circle*{2}}
\put(130,40){\circle*{2}}
\put(150,40){\circle*{2}}

\put(10,20){\line(0,1){20}}
\put(30,20){\line(0,1){20}}
\put(70,20){\line(1,1){20}}
\put(90,20){\line(-1,1){8}}
\put(70,40){\line(1,-1){8}}

\put(130,20){\line(0,1){20}}
\put(150,20){\line(0,1){20}}

\put(45,25){$\cdots$}
\put(105,25){$\cdots$}

\put(8,5){\scriptsize $1$ }
\put(28,5){\scriptsize $2$ }
\put(68,5){\scriptsize $j$ }
\put(86,5){\scriptsize $j+1$ }
\put(122,5){\scriptsize $2n-1$ }
\put(150,5){\scriptsize $2n$ }

\end{picture}
\end{center}

\vspace{-1em}

\noindent
then each $\ga_{2i-1}$ is just an elementary
braid on the doubled strands:

\vspace{-2em}

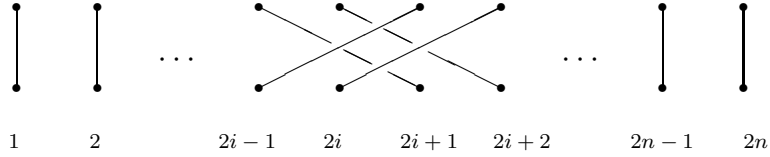
\setlength{\unitlength}{0.13em}
\begin{center}
\begin{picture}(240,60)

\thinlines

\put(10,20){\circle*{2}}
\put(30,20){\circle*{2}}
\put(70,20){\circle*{2}}
\put(90,20){\circle*{2}}
\put(110,20){\circle*{2}}
\put(130,20){\circle*{2}}
\put(170,20){\circle*{2}}
\put(190,20){\circle*{2}}

\put(10,40){\circle*{2}}
\put(30,40){\circle*{2}}
\put(70,40){\circle*{2}}
\put(90,40){\circle*{2}}
\put(110,40){\circle*{2}}
\put(130,40){\circle*{2}}
\put(170,40){\circle*{2}}
\put(190,40){\circle*{2}}

\put(10,20){\line(0,1){20}}
\put(30,20){\line(0,1){20}}
\put(70,20){\line(2,1){40}}
\put(90,20){\line(2,1){40}}
\put(70,40){\line(2,-1){18}}
\put(90,40){\line(2,-1){8}}

\put(110,20){\line(-2,1){8}}
\put(130,20){\line(-2,1){18}}

\put(98,26){\line(-2,1){7}}
\put(108,31){\line(-2,1){7}}

\put(170,20){\line(0,1){20}}
\put(190,20){\line(0,1){20}}

\put(45,25){$\cdots$}
\put(145,25){$\cdots$}

\put(8,5){\scriptsize $1$ }
\put(28,5){\scriptsize $2$ }
\put(60,5){\scriptsize $2i-1$ }
\put(86,5){\scriptsize $2i$ }
\put(105,5){\scriptsize $2i+1$ }
\put(128,5){\scriptsize $2i+2$ }
\put(162,5){\scriptsize $2n-1$ }
\put(190,5){\scriptsize $2n$ }

\end{picture}
\end{center}

\vspace{-1em}

\noindent
For each odd $i$ the elements \eqref{diabl} generate a subalgebra
of $\U'_q(\spa_{2n})$ isomorphic to $\U'_q(\spa_{2})$.
The next proposition shows that the elements $\ga_i$ permute
these subalgebras.

\bpr\label{prop:bgodd}
The images of the elements \eqref{diabl}
under the action of the automorphisms
$\ga_1,\ga_3,\dots,\ga_{2n-3}$ belong to
$\U'_q(\spa_{2n})$.
\epr

\bpf
This is verified with the use of \eqref{sijtij}.
For any odd $j$
the elements \eqref{diabl} with $i\ne j,j+2$ are fixed by the action
of $\ga_j$, while
\ben
\ga_j: s^{}_{jj}\mapsto s^{}_{j+2,j+2},
\qquad s^{}_{j+1,j+1}\mapsto s^{}_{j+3,j+3},\qquad
s^{}_{j,j+1}\mapsto s^{}_{j+2,j+3}
\een
and
\ben
\ga_j: s^{}_{j+2,j+2}\mapsto s^{}_{jj},
\qquad s^{}_{j+3,j+3}\mapsto s^{}_{j+1,j+1},\qquad
s^{}_{j+2,j+3}\mapsto s^{}_{j,j+1}.
\een
This follows from the formulas for the action of the $\be_i$ on
$\U_q(\gl_{2n})$ which imply, for instance, relations
of the type
\ben
\be_{j}\be_{j+1}:t_{j+1,j}\mapsto t_{j+2,j+1}.
\een
Since $\ga_j=\be_{j+1}\be_{j+2}\be_{j}\be_{j+1}$, this gives
$\ga_j:t_{j+1,j}\mapsto t_{j+3,j+2}$. The images of the remaining
elements of the form $t_{jj},\bar t_{j,j+1},t_{j+1,j+1}$ are calculated
in a similar way which gives the desired formulas.
\epf

It can be shown that Proposition~\ref{prop:bgodd} is not
extended to the remaining generators
\eqref{undial} of the algebra $\U'_q(\spa_{2n})$.
Observe that the elements $\be_i$ and $\ga_i$ of $B_{2n}$
with odd $i$
satisfy the relations
\ben
\ga_i^{-1}\be^{}_j\ts\ga^{}_i=\be^{}_j\qquad\text{if}\quad j\ne i,i+2
\een
while
\ben
\ga_i^{-1}\be^{}_i\ts\ga^{}_i=\be^{}_{i+2}\Fand
\ga_i^{-1}\be^{}_{i+2}\ts\ga^{}_i=\be^{}_{i}.
\een
The elements $\be_i$ generate a subgroup of $B_{2n}$
isomorphic to $\ZZ^n$. We shall identify $\ZZ^n$ with this subgroup.
These observations suggest the following definition. Consider the braid
group $B_n$ with generators $\ga'_1,\ga'_3,\dots,\ga'_{2n-3}$
and the usual defining relations
\ben
\ga'_i\tss \ga'_{i+2}\tss \ga'_i=\ga'_{i+2}\tss \ga'_i\tss \ga'_{i+2},
\qquad i=1,3,\dots,2n-5
\een
and
\ben
\ga'_i\tss \ga'_j=\ga'_j\tss \ga'_i,\qquad |i-j|>2.
\een
Define the group $\Gamma_n$ as the semidirect product
$\Gamma_n=B_n\ltimes \ZZ^n$ where the action of $B_n$ on $\ZZ^n$
is defined by
\ben
\be_j^{\ga'_i}=\be^{}_j\qquad\text{if}\quad j\ne i,i+2
\een
while
\ben
\be_i^{\ga'_i}=\be^{}_{i+2}\Fand
\be_{i+2}^{\ga'_i}=\be^{}_{i}.
\een
Note that the Weyl group $W(C_n)=\Sym_n\ltimes \ZZ_2^n$ of type $C_n$
may be regarded as a classical counterpart of  $\Gamma_n$.

\bcj\label{conj:act}
There exists an action of the group $\Gamma_n$ on the algebra
$\U'_q(\spa_{2n})$ by automorphisms which corresponds to
the action of $W(C_n)$ on $\U(\spa_{2n})$.
\ecj

Our final theorem shows that the conjecture holds for $n=2$.

\bth\label{thm:netwo}
Let the generators $\be_1$ and $\be_3$ of
the group $\Gamma_2$ act on $\U'_q(\spa_{4})$
as in Proposition~\ref{prop:bodd} and let the generator
$\ga'_1$ act on the elements \eqref{diabl}
with $i=1,3$ as $\ga_1$.
Then together with the assignment
\ben
\ga'_1:s_{32}\mapsto s_{41},\qquad s_{41}\mapsto s_{32},
\qquad s_{31}\mapsto s_{31},\qquad s_{42}\mapsto s_{42}
\een
this defines an action of $\Gamma_2$ on $\U'_q(\spa_{4})$
by automorphisms.
\eth

\bpf
It is easy to verify that $\ga'_1$ respects the defining relations
of $\U'_q(\spa_{4})$. For instance, the following relations
are clearly respected by $\ga'_1$
\ben
\bal
s_{33}\ts s_{32}&=s_{32}\ts s_{33},
\qquad
s_{11}\ts s_{32}=s_{32}\ts s_{11}+(q^{-1}-q)\ts s_{12}\ts s_{31}\\
s_{31}\ts s_{32}&=q^{-1}\ts s_{32}\ts s_{31}+(q-\qin)
(q^{-1}\ts s_{21}\ts s_{33}-s_{12}\ts s_{33})
\eal
\een
and
\ben
\bal
s_{11}\ts s_{41}&=s_{41}\ts s_{11},
\qquad
s_{33}\ts s_{41}=s_{41}\ts s_{33}+(q^{-1}-q)\ts s_{34}\ts s_{31}\\
s_{31}\ts s_{41}&=q^{-1}\ts s_{41}\ts s_{31}+(q-\qin)
(q^{-1}\ts s_{43}\ts s_{11}-s_{34}\ts s_{11})
\eal
\een
together with
\ben
s_{32}\ts s_{41}=s_{41}\ts s_{32}+(q-\qin)(s_{12}\ts s_{43}-s_{34}\ts s_{21}),
\een
and this holds for
the remaining relations as well.
The defining relations of the group $\Gamma_2$ are also easily verified.
\epf

\section*{Appendix}\label{sec:app}
\setcounter{section}{5}
\setcounter{equation}{0}

Here we prove that the Sklyanin determinant
$\sdet S(u)$ is central in
the extended algebra $\hat\U'_q(\spa_{N})$ with $N=2n$;
see the proof of Theorem~\ref{thm:casimsym}.
We need to introduce some
more notation. Following \cite{mrs:cs}, introduce
the trigonometric $R$-matrix
\beql{trRm}
\bal
R(u,v)={}&(u-v)\sum_{i\ne j}E_{ii}\ot E_{jj}+(\qin u-q\tss v)
\sum_{i}E_{ii}\ot E_{ii} \\
{}+ {}&(\qin-q)\tss u\tss\sum_{i> j}E_{ij}\ot
E_{ji}+ (\qin-q)\tss v\tss\sum_{i< j}E_{ij}\ot E_{ji}
\eal
\eeq
and a rational function in independent variables
$u_1,\dots,u_r,q$ valued in $(\End\CC^N)^{\ot\tss r}$ by
\beql{Rlong}
R(u_1,\dots,u_r)=\prod_{i<j}R_{ij}(u_i,u_j),
\eeq
where
the product is taken in the lexicographical order on the pairs $(i,j)$.
We have the following relation in
the algebra $\hat\U'_q(\spa_{N})\ot(\End\CC^N)^{\ot\tss r}$,
\begin{multline}\label{fundamtw}
R(u_1,\dots,u_r)\ts S_1(u_1)R_{12}^{\tss t}\cdots R_{1r}^{\tss t}   S_2(u_2)
R_{23}^{\tss t}\cdots R_{2r}^{\tss t}   S_3(u_3) \cdots R_{r-1,r}^{\tss t}
S_r(u_r)=\\
S_r(u_r)R_{r-1,r}^{\tss t}\cdots S_3(u_3)   R_{2r}^{\tss t} \cdots R_{23}^{\tss t}
S_2(u_2)R_{1r}^{\tss t} \cdots R_{12}^{\tss t} S_1(u_1)\ts R(u_1,\dots,u_r);
\end{multline}
see \cite{mrs:cs},
where $R_{ij}^{\tss t}=R_{ij}^{\tss t}(u_i^{-1},u_j)$ with $R^{\tss t}(u,v)$
defined in \eqref{rtuv}. Now take $r=N+1$ and label
the copies of $\End\CC^N$ in the tensor product
$\hat\U'_q(\spa_{N})\ot(\End\CC^N)^{\ot\tss (N+1)}$ with the indices
$0,1,\dots,N$. Furthermore, specialize
the parameters $u_i$ in \eqref{fundamtw} as follows:
\ben
u_0=v, \qquad u_i=q^{-2i+2}u\quad\text{for}\ \ i=1,\dots,N.
\een
Then by \cite[Proposition~4.1]{mrs:cs},
the element \eqref{Rlong} will take the form
\ben
R(v,u,\dots,q^{-2N+2}u)=\al(u)\ts\prod_{i=1,\dots,N}^{\longrightarrow}
R_{0i}(v,q^{-2i+2}u)\ts A^q_N,
\een
where
\ben
\al(u)=u^{N(N-1)/2}\ts
\prod_{1\leqslant i<j\leqslant N}(q^{-2i+2}-q^{-2j+2}).
\een
We shall now be verifying that
\beql{ancoll}
\prod_{i=1,\dots,N}^{\longrightarrow}
R_{0i}(v,q^{-2i+2}u)\ts A^q_N=\de(u,v)\ts A^q_N
\eeq
where
\ben
\de(u,v)=(\qin v-q\ts u)\prod_{i=1}^{N-1}(v-q^{-2i}u).
\een
The $R$-matrix $R(u,v)$ satisfies the Yang--Baxter equation
\ben
R_{12}(u,v)  R_{13}(u,w)R_{23}(v,w) =  R_{23}(v,w) R_{13}(u,w) R_{12}(u,v).
\een
Using this relation repeatedly, we derive the identity
\ben
R(u_1,\dots,u_r)=\prod_{i<j}R_{ij}(u_i,u_j),
\een
where the product is taken in the order opposite to the
lexicographical order on the pairs $(i,j)$.
Taking here $r=N+1$ and specializing the variables $u_i$ as above,
we arrive at
\beql{symma}
\prod_{i=1,\dots,N}^{\longrightarrow}
R_{0i}(v,q^{-2i+2}u)\ts A^q_N=A^q_N\ts \prod_{i=1,\dots,N}^{\longleftarrow}
R_{0i}(v,q^{-2i+2}u).
\eeq
Hence, for the proof of \eqref{ancoll}, it now suffices
to compare the images of the operators on both sides
at the basis vectors of the form $v_k=e_k\ot e_{i_1}\ot\cdots\ot\ts e_{i_N}$
with $k=1,\dots,N$,
where the $e_i$ denote the canonical basis vectors of $\CC^N$
and $\{i_1,\dots,i_N\}$ is a fixed permutation of $\{1,\dots,N\}$.
Our next observation is the fact that for any $i,j\in\{1,\dots,N\}$
the expression
$R(u,v)(e_i\ot e_j)$ is a linear combination
of $e_i\ot e_j$ and $e_j\ot e_i$. This implies that for each $k$,
\beql{aqn}
A^q_N\ts \prod_{i=1,\dots,N}^{\longleftarrow}
R_{0i}(v,q^{-2i+2}u)\ts v_k=\de_k(u,v)\ts A^q_N\ts v_k,
\eeq
for some scalar function $\de_k(u,v)$ which is independent of
the permutation $\{i_1,\dots,i_N\}$.
It remains to show that
$\de_k(u,v)=\de(u,v)$ for all $k$. However, this is immediate from
\eqref{aqn} if for a given $k$ we choose a permutation $\{i_1,\dots,i_N\}$
with $i_1=k$, thus completing the proof of \eqref{ancoll}.

Now apply the transposition $t$ on the $0$-th copy of $\End\CC^N$
and combine \eqref{ancoll} and \eqref{symma}
to derive another identity
\ben
A^q_N\ts\prod_{i=1,\dots,N}^{\longrightarrow}
R_{0i}^{\tss t}(v,q^{-2i+2}u)
=\prod_{i=1,\dots,N}^{\longleftarrow}
R_{0i}^{\tss t}(v,q^{-2i+2}u)\ts A^q_N=\de(u,v)\ts A^q_N.
\een
Thus, \eqref{fundamtw} becomes
\ben
\de(u,v)\ts \de(u,v^{-1})\ts A^q_N\ts S_0(v)\ts \sdet S(u)=
\de(u,v)\ts \de(u,v^{-1}) \ts A^q_N\ts \sdet S(u)\ts S_0(v),
\een
proving that $\sdet S(u)$ lies in the center of $\hat\U'_q(\spa_{N})$.

As a final remark, note that the above argument applies
to more general matrices $S(u)$.
The only property of $S(u)$ used above
is the fact that $S(u)$ satisfies the reflection equation
\beql{refle}
R(u,v)\ts S_1(u)\ts R^{\ts t}(u^{-1},v)\ts S_2(v)=S_2(v)
\ts R^{\ts t}(u^{-1},v)\ts S_1(u)\ts R(u,v).
\eeq
This implies that \eqref{sdetmatr} equals $A^q_N\ts\sdet S(u)$
for some formal series $\sdet S(u)$ called
the Sklyanin determinant. Then $\sdet S(u)$ is central
in the algebra with the defining relations \eqref{refle}.
In particular, this applies to the (extended) twisted $q$-Yangians
associated with the orthogonal and symplectic Lie algebras;
see \cite{mrs:cs}.

\end{document}